\documentclass[english,11pt]{amsart}
\usepackage{amsopn,amsthm,amsfonts,amsmath,amssymb,amscd}
\usepackage{babel}
\usepackage{nicematrix}
\usepackage{mathtools}
\usepackage{tikz}

\newtheorem{Theorem}{Theorem}[section]
\newtheorem{Lemma}[Theorem]{Lemma}

\newtheorem{Example}[Theorem]{Example}
\newtheorem*{Remark}{Remark}

\newenvironment{Proof}{{\bf Proof.}}{\hfill $\blacksquare$}
\newenvironment{Proof*}{{\it Proof.}}

\newcommand{\CC}{\mathbb{C}}
\newcommand{\RR}{\mathbb{R}}
\newcommand{\ZZ}{\mathbb{Z}}

\newcommand{\diag}{{\rm diag}}
\newcommand{\antidiag}{{\rm antidiag}}

\newcommand{\rk}{\mathrm{rk}}

\newcommand{\ch}{{\rm char}}

\begin{document}

\title{Eigenvalues of the product matrices of finite commutative rings}

\author{David Dol\v zan}

\address{D.~Dol\v zan:~Department of Mathematics, Faculty of Mathematics
and Physics, University of Ljubljana, Jadranska 19, SI-1000 Ljubljana, Slovenia, and Institute of Mathematics, Physics and Mechanics, Jadranska 19, SI-1000 Ljubljana, Slovenia}

\email{david.dolzan@fmf.uni-lj.si}

\thanks{The author acknowledges the financial support from the Slovenian Research Agency  (research core funding No. P1-0222)}

\begin{abstract} 
The product matrix of a finite commutative ring $R=\{x_1,x_2,\ldots,x_n\}$ and an element $u \in R$ is the matrix $A_u(R)=[a_{ij}]$, where $a_{ij}=1$ if $x_ix_j=u$, and $a_{ij}=0$ otherwise. This provides a natural extension of the concept of the adjacency matrix of the zero-divisor graph of a ring, which has been studied extensively. In this paper, we find the characteristic polynomial of $A_u(R)$ for a local ring $R$ of odd order and a unit $u$. By studying the structure of a finite local ring, we find the characteristic polynomial of $A_u(R)$ for a local ring $R$ and any $u \in R$ in two cases: when the Jacobson radical of $R$ has either the maximal or the minimal possible index of nilpotency.
\end{abstract}

\maketitle 

 \section{Introduction}

\bigskip

The structure of the zero-divisors has always played an important role in ring theory and this topic has been heavily studied in recent years, providing many classical results in the study of finite rings (see, for example, \cite{ganesan, ganesan1, koh1, koh}). Later, many authors concerned themselves with probabilities of zero multiplication (sometimes also called the nullity degree) of finite rings (see, for example, \cite{dolzan, esmkhani, esmkhani1, salih, shumyatsky}). This topic has been extended to studying the probabilities that a product of two elements is equal to a certain fixed element of the ring (instead of zero) --- see, for example, \cite{dolzan1, rehman, rehman0}.

A natural way to study the zero-divisors in a ring is to associate to the ring its zero-divisor graph (i.e., a graph with vertices corresponding to all non-zero zero-divisors and edges between any two vertices that multiply to zero). Some authors also study the extended zero-divisor graph (where the set of vertices consists of all elements of the ring). It turns out that one of the most fruitful areas of research in this field is studying the spectra of different matrices that correspond to the zero-divisor graphs.  For example, authors in \cite{bajaj, barati, chatto, monius, pirzada1, pirzada2, pirzada, rather, ratta} (among many others) studied the eigenvalues of the adjacency matrix (and also the Laplacian matrix) of the zero-divisor graph of the ring of integers modulo $n$.

In this paper, we extend the above in a natural way --- namely, we study the eigenvalues of the product matrix corresponding to an arbitrary element in a finite ring. To be precise, if $R$ is a finite commutative ring with elements $\{x_1, x_2, \ldots, x_n\}$ and $u \in R$, we define the $n$-by-$n$ matrix $A_u(R)=[a_{ij}]$ with $a_{ij}=\begin{cases}
1; \text { if } x_ix_j=u, \\
0; \text { otherwise},    
\end{cases}$ and we call this matrix the \emph{product matrix of the element} $u \in R$. Obviously, if $u=0$, the product matrix is exactly the adjacency matrix of the extended zero-divisor graph. We find the characteristic polynomial (and thus the eigenvalues) of this matrix for some special cases of local rings, thereby generalizing many of the previously known results in two ways --- firstly, by studying the product matrices instead of the adjacency matrices of the zero-divisor graphs, and secondly by considering more general classes of local rings rather than just the rings of integers modulo $n$.

The structure of this paper is as follows. In the next section, we gather the necessary definitions and preliminary results that we shall use throughout the paper. In Section 3, we find the characteristic polynomial of the product matrix $A_u(R)$ for an arbitrary commutative local ring $R$ of odd order in the case when $u$ is a unit in $R$ (see Theorem \ref{xjeenota}).
In Section 4, we study the local commutative rings such that the Jacobson radical has the maximal possible index of nilpotency. We find the characteristic polynomial of  $A_0(R)$ (see Theorem \ref{xjenula}), thereby proving similar  results (but in a more general setting and using different methods) for the eigenvalues of the zero-divisor graphs as in \cite[Corollary 2.8 and Corollary 2.13]{pirzada}, where the eigenvalues of the zero-divisor graphs (without loops) of the ring $\ZZ_{p^nq^m}$ have been studied. We also find the characteristic polynomial of $A_u(R)$ for a commutative local ring $R$ of even order in the case when $u$ is a unit in $R$ (see Theorem \ref{xjeenotaeven}) and in the case when $u$ is a non-zero zero-divisor (see Theorem \ref{xjedn}). In the final section, we then find the characteristic polynomial of $A_u(R)$ in the case of a local commutative ring such that its Jacobson radical has the minimal possible index of nilpotency (see Theorem \ref{jkvadratnula}).

\bigskip

 \section{Definitions and preliminaries}

Throughout the paper, $R$ will denote a finite commutative ring with unity $1$, its group of units will be denoted by $R^*$, its set of zero-divisors by $Z(R)$ and its Jacobson radical by $J(R)$. When the ring in question will be beyond doubt, we will write simply $J$ instead of $J(R)$. We shall often denote the equivalence class of $x$ in $R/J$ by $\overline{x}=x+J$. 



We will denote by $M_n(F)$ the ring of $n$-by-$n$ matrices over the field $F$ and we will denote the rank of $A \in M_n(F)$ by $\rk(A)$.

For a matrix $A \in M_n(F)$ we will denote the $(i,j)$-th entry of $A$ by $A_{i,j}$.
The identity matrix in $M_n(F)$ will be denoted by $I_n$, the zero matrix in $M_n(F)$ will be denoted by $0_n$ and the $i$-by-$j$ zero matrix will be denoted by $0_{i,j}$.  Furthermore, let $F_{i,j}$ denote the $i$-by-$j$ matrix with all elements equal to $1$, and we will also use the notation $F_i$ for $F_{i,i}$.
For matrices $A \in M_{k}(F)$ and $B \in M_l(F)$ we shall denote their direct sum by $A \oplus B \in M_{k+l}(F)$, that is $A \oplus B =  \begin{pNiceMatrix}
    A & 0_{k,l} \\

0_{l,k} &  {B} \end{pNiceMatrix}$, and $A \otimes B$ will denote the Kronecker product of the matrices $A$ and $B$. Finally, the diagonal matrix with elements $a_1,a_2,\ldots,a_n$ along the diagonal will be denoted by $\diag(a_1,a_2,\ldots,a_n)$ and similarly, $\antidiag(a_1,a_2,\ldots,a_n)$ will denote the matrix $A$ with all its elements equal to zero, except for the elements $A_{1,n}=a_1, A_{2,n-1}=a_2, \ldots, A_{n,1}=a_n$.

For a matrix $A \in M_n(F)$, its characteristic polynomial is $p(\lambda)=\det(A-\lambda I)$.


The next lemma will be useful.

\begin{Lemma}
\label{inver}
  Let $R$ be a (finite) commutative ring, $u \in R^*$ and $j \in J$. Then $u+j \in R^*$.
\end{Lemma}
\begin{Proof}
 This follows directly from the definition of the Jacobson radical.
\end{Proof}

\bigskip

We say that a (finite) ring $R$ is a \emph{local ring} if $R$ has a unique maximal left ideal. For finite rings, this is equivalent (see, for example, \cite[Theorem V.1]{mcdonald}) to the fact that $R/J$ is a field and also to the fact that the non-units of $R$ form an additive abelian group.

\bigskip

Obviously, an example of a finite local ring is the ring $\ZZ_{p^n}$ for any prime number $p$ and any integer $n$. We shall use the following throughout the paper.

\begin{Theorem}\cite[Theorem 2]{raghavendran}
  Let $R$ be a finite local ring. Then there exists a prime number $p$ and integers $n,r$ such that $|R|=p^{nr}$, $|J|=p^{(n-1)r}$ and $R/J$ is a field with $p^r$ elements.   Furthermore, $\ch(R)=p^k$ for some $1 \leq k \leq n$ and $J^n=0$.
\end{Theorem}

\bigskip

 \section{Local rings of odd order}

\bigskip

As we have noted in the introduction, if $u = 0$, then the product matrix $A_0(R)$ is exactly the adjacency matrix of the extended zero-divisor graph of the ring $R$. However, we shall deal with this case in Section 4 (see Theorem \ref{xjenula}).
In this section, we shall start with a simpler case --- we shall study the product matrix $A_u(R)$ in the case when $u \in R^*$ and $R$ is a commutative local ring of odd order.
We have the following theorem.

\begin{Theorem}
  \label{xjeenota}
   Let $R$ be a finite commutative local ring of order $q^n$ such that $R/J$ is a field of odd order $q$. For $u \in R^*$, let $p(\lambda)$ denote the characteristic polynomial of the matrix $A_u(R)$. If $u=s^2$ for some $s \in R$, then $p(\lambda)=-\lambda^{q^{n-1}}(\lambda-1)^{\frac{q^n-q^{n-1}+2}{2}}(\lambda+1)^{\frac{q^n-q^{n-1}-2}{2}}$. For all other $u \in R^*$, we have $p(\lambda)=-\lambda^{q^{n-1}}(\lambda-1)^{\frac{q^n-q^{n-1}}{2}}(\lambda+1)^{\frac{q^n-q^{n-1}}{2}}$.
   Furthermore, there exist exactly $\frac{q^n-q^{n-1}}{2}$ elements $u \in R^*$ such that $u=s^2$ for some $s \in R$.
\end{Theorem}
\begin{Proof}
   Obviously, for $s \in J$ there exists no $t \in R$ such that $st=u \in R^*$.
   On the other hand, if $s \in R^*$, we have $s(s^{-1}u)=u$. Observe that for any $s \in R^*$, the fact that $s'(s^{-1}u) = u$ for some $s' \in R^*$ implies that $s'=s$, so for every $s \in R^*$ there exists exactly one $t \in R$ such that $st=u$. Note also that $s^2=u$ for $u \in R^*$ implies that $s \in R^*$, so that $s^{-1}u=s$ if and only if $s^2=u$. 
   Suppose now that $s^2=v^2$ for some $s, v \in R^*$. Since $R/J$ is a field of odd order, this implies $s=v+j$ or $s=-v+j$ for some $j \in J$.
   Say, $s=v+j$ (we treat the other case similarly). Then $s^2=v^2+2vj+j^2=v^2$, so $(2v+j)j=0$. Since $q$ is odd, $2 \notin J$, so $2$ is an invertible element in $R$. Then by Lemma \ref{inver}, $2v+j$ is also invertible, which implies that $j=0$ and $s=v$.
Therefore, we have proved that for any $s \in R^*$ there exist exactly two $v \in R^*$ such that $s^2=v^2$. The mapping $\alpha : R^* \rightarrow R^*$ defined by $\alpha(x)=x^2$, thus has an image of order $\frac{|R^*|}{2}=\frac{q^n-q^{n-1}}{2}$.
   Now, denote $B_k=\begin{pNiceMatrix}
    0 & 1 \\

1 &  0 \end{pNiceMatrix} \oplus \begin{pNiceMatrix}
    0 & 1 \\

1 &  0 \end{pNiceMatrix} \oplus \ldots \oplus \begin{pNiceMatrix}
    0 & 1 \\

1 &  0 \end{pNiceMatrix} \in M_{2k}(\RR)$ ($k$ summands).
   This means that in the case when $u$ is a square, we have (with the following reordering of elements in $R$: $q^{n-1}$ elements from $J$, $2$ elements $s \in R \setminus J$ such that $s^2=u$, and then the other $q^n-q^{n-1}-2$ invertible elements from $R$ in pairs $s$, $s^{-1}u$),
    $A_u(R)=0_{q^{n-1}} \oplus I_{2} \oplus B_\frac{q^n-q^{n-1}-2}{2}$. Since the spectrum of the matrix $\begin{pNiceMatrix}
    0 & 1 \\

1 &  0 \end{pNiceMatrix}$ is $\{-1,1\}$, we have $p(\lambda)=-\lambda^{q^{n-1}}(\lambda-1)^{\frac{q^n-q^{n-1}+2}{2}}(\lambda+1)^{\frac{q^n-q^{n-1}-2}{2}}$. In the case when $u$ is not a square, we have (with a similar reordering),
$A_u(R)=0_{q^{n-1}} \oplus B_\frac{q^n-q^{n-1}}{2}$, therefore $p(\lambda)=-\lambda^{q^{n-1}}(\lambda-1)^{\frac{q^n-q^{n-1}}{2}}(\lambda+1)^{\frac{q^n-q^{n-1}}{2}}$.
\end{Proof}

\bigskip

 \section{Local rings with the Jacobson radical of maximal nilpotency index}

\bigskip

In this section, we will concern ourselves with the following setup: $R$ will be a finite commutative ring of order $q^n$, $R/J$ a field of order $q$ and $J^{n-1} \neq 0$. Observe that this class of rings contains all rings $\ZZ_{p^n}$ for any prime $p$ and any integer $n$, and thereby the results of this section will improve upon some results regarding the spectrum of the adjacency matrix of the zero-divisor graph of  $\ZZ_{p^n}$.

\bigskip

%


\begin{Theorem}
  \label{xjenula}
   Let $R$ be a finite commutative local ring of order $q^n$, $R/J$ a field of order $q$ and $J^{n-1} \neq 0$. Then the characteristic polynomial of the matrix $A_0(R)$ is equal to $p(\lambda)=(-1)^{q-\frac{(n+1)(n+2)}{2}}(q-1)^nq^\frac{n(n-1)}{2}\lambda^{q^n-(n+1)}\det(B-\lambda \cdot\antidiag(1,1/\alpha_1,1/\alpha_2,\ldots,1/\alpha_{n}))$, where $B \in M_{n+1}(\RR)$ is the upper-triangular matrix with all its elements equal to $1$ and $\alpha_i=q^i-q^{i-1}$ for $i=1,2,\ldots,n$.
\end{Theorem}
\begin{Proof}
    We order the elements of $R$ in the following way: units, elements from $J \setminus J^2$, elements from $J^2 \setminus J^3$, \ldots, elements from $J^{n-1} \setminus \{0\}$, $0$. Observe that in this ordering, we have $$A_0(R)= \begin{pNiceMatrix}[margin,columns-width=auto,nullify-dots]
    0_{\alpha_n,\alpha_n} & 0_{\alpha_n,\alpha_{n-1}} & 0_{\alpha_n,\alpha_{n-2}} & \ldots & 0_{\alpha_n,\alpha_1} & F_{\alpha_n,1} \\
    0_{\alpha_{n-1},\alpha_{n}} & 0_{\alpha_{n-1},\alpha_{n-1}} & 0_{\alpha_{n-1},\alpha_{n-2}}  & \ldots & F_{\alpha_{n-1},\alpha_{1}} & F_{\alpha_{n-1},1} \\
    0_{\alpha_{n-2},\alpha_{n}} & 0_{\alpha_{n-2},\alpha_{n-1}} &   & \iddots & F_{\alpha_{n-2},\alpha_{1}} & F_{\alpha_{n-2},1} \\
    \\
\vdots & & \iddots & \iddots & \vdots &  \vdots \\
\\
0_{\alpha_{1},\alpha_{n}} & F_{\alpha_{1},\alpha_{n-1}} & F_{\alpha_{1},\alpha_{n-2}}  & \ldots & F_{\alpha_{1},\alpha_{1}} & F_{\alpha_{1},1} \\
F_{1,\alpha_{n}} & F_{1,\alpha_{n-1}} & F_{1,\alpha_{n-2}}  & \ldots & F_{1,\alpha_{1}} & 1  \\
\end{pNiceMatrix},$$ where $\alpha_i=q^i-q^{i-1}$ for all $i=1,2,\ldots,n$. 
Note that $\rk(A_0(R))=n+1$, and since $A_0(R)$ is diagonalizable ($R$ is commutative, so $A_0(R)$ is symmetric), $0$ is an eigenvalue of $A_0(R)$ with the algebraic multiplicity  equal to $q^n-(n+1)$.
For every $i=1,2,\ldots,n$ denote by $v_i \in \RR^{q^n}$ the vector with its last $1+\alpha_1+\alpha_2+\ldots+\alpha_i$ components equal to $1$ and all the other components equal to $0$. Furthermore, let $v_0 \in \RR^{q^n}$ be the vector with all its components equal to $0$ except for the last component which is equal to $1$. Since for any $v \in \CC^{q^n}$, the equation $A_0(R)x=v$ is solvable if and only if $v$ lies in the column space of $A_0(R)$, we see that every eigenvector of $A_0(R)$ that corresponds to a non-zero eigenvalue is a linear combination of $v_0,v_1, \ldots, v_n$. Furthermore, it is easy to see that $$A_0(R)v_i=v_n+\sum_{j=n-i}^{n-1}{\alpha_{n-j}v_j}$$ for every $i=0,1,\ldots,n$. Assume now that $w=\sum_{i=0}^n{\beta_iv_i}$ for some $\beta_i \in \CC$ is an eigenvector of $A_0(R)$ corresponding to a non-zero eigenvalue $\lambda$. Then $A_0(R)w=\sum_{i=0}^n{\beta_iA_0(R)(v_i)}=\sum_{i=0}^n{\beta_i}v_n+\sum_{j=n-i}^{n-1}{\alpha_{n-j}v_j}=\sum_{i=0}^n{\lambda\beta_iv_i}$. This leads to the following equations:
\begin{equation}
\label{sistem}
\lambda\beta_n=\sum_{i=0}^n{\beta_i},  \,\,\, \lambda\beta_{n-j}=\sum_{i=j}^n{\beta_i\alpha_j} \text{ for } j=1,2,\ldots,n. 
\end{equation}
Denote $z=(\beta_n,\beta_{n-1},\ldots,\beta_0)^T \in \CC^{n+1}$,  $B'=\begin{pNiceMatrix}[margin,columns-width=auto,nullify-dots]
    1 & 1 & 1 & \ldots & 1 & 1 \\
    0 & \alpha_1 & \alpha_1 & \ldots & \alpha_1 & \alpha_1 \\
    0 & 0 & \alpha_2 & \ldots & \alpha_2 & \alpha_2 \\
\vdots & & \ddots & \ddots & \vdots &  \vdots \\
    0 & 0 & 0 & \ddots & \alpha_{n-1} & \alpha_{n-1} \\
    0 & 0 & 0 & \ldots & 0 & \alpha_{n} \\
    \end{pNiceMatrix}$ and $C'=\begin{pNiceMatrix}[margin,columns-width=auto,nullify-dots]
    0 & 0 &  \ldots & 0 & \lambda \\
    0 & 0 &  & \lambda & 0 \\
\vdots & & \iddots & & \vdots \\
0 & \lambda &   & 0 & 0 \\
    \lambda &  0 & \ldots & 0 & 0 \\
    \end{pNiceMatrix} \in M_{n+1}(\CC)$ and observe that equations (\ref{sistem}) transform into the equation $(B'-C')z=0$. Since we are searching for non-trivial solutions to equations (\ref{sistem}), we have $\det(B'-C')=0$. Now, if we divide the $i$-th row of $B'-C'$ by $\alpha_i$ for $i=1,2,\ldots,n$, we get $\det(B- \lambda \cdot \antidiag(1,1/\alpha_1,1/\alpha_2,\ldots,1/\alpha_{n}))=0$.  If we denote $D= \antidiag(1,1/\alpha_1,1/\alpha_2,\ldots,1/\alpha_{n})$, we have     $\det(B- \lambda  D)=\det(D)\det(D^{-1}B-\lambda I_{n+1})$, so the leading coefficient of the polynomial $\det(B- \lambda D)$ equals $(-1)^{\frac{(n+1)(n+2)}{2}}\prod_{i=1}^n\frac{1}{\alpha_i}$. Since $\prod_{i=1}^n{\alpha_i}=\prod_{i=1}^n{q^{i-1}(q-1)}=(q-1)^nq^\frac{n(n-1)}{2}$, 
    the theorem is now proven.
\end{Proof}

\begin{Example}
Let us illustrate this theorem in a simple case of $R=\ZZ_{q^2}$ for a prime $q$.
Theorem \ref{xjenula} gives us that the characteristic polynomial of the matrix $A_0(R)$ is equal to $p(\lambda)=(-1)^{q}\lambda^{q^2-3}(\lambda^3-q\lambda^2-q(q-1)\lambda+q(q-1)^2)$.
\end{Example}

\begin{Remark}
    Observe that Theorem \ref{xjenula} also describes the characteristic polynomial of the adjacency matrix of the zero-divisor graph of $R$, if we allow the zero-divisor graph to have loops.  Namely, if we denote the adjacency matrix of the zero-divisor graph (with loops) of $R$ by $A$, we can see that $A$ is the principal submatrix of $A_0(R)$, obtained by striking out the rows and columns that correspond to the elements in $R^* \cup \{0\}$. But then $A=\begin{pNiceMatrix}[margin,columns-width=auto,nullify-dots]
    0_{\alpha_{n-1},\alpha_{n-1}} & 0_{\alpha_{n-1},\alpha_{n-2}} & 0_{\alpha_{n-1},\alpha_{n-3}} & \ldots & 0_{\alpha_{n-1},\alpha_2} & F_{\alpha_{n-1},\alpha_1} \\
    0_{\alpha_{n-2},\alpha_{n-1}} & 0_{\alpha_{n-2},\alpha_{n-2}} & 0_{\alpha_{n-1},\alpha_{n-3}}  & \ldots & F_{\alpha_{n-2},\alpha_{2}} & F_{\alpha_{n-2},\alpha_1} \\
    0_{\alpha_{n-3},\alpha_{n-1}} & 0_{\alpha_{n-3},\alpha_{n-2}} &   & \iddots & F_{\alpha_{n-3},\alpha_{2}} & F_{\alpha_{n-3},\alpha_1} \\
    \\
\vdots & & \iddots & \iddots & \vdots &  \vdots \\
\\
0_{\alpha_{2},\alpha_{n-1}} & F_{\alpha_{2},\alpha_{n-2}} & F_{\alpha_{2},\alpha_{n-3}}  & \ldots & F_{\alpha_{2},\alpha_{2}} & F_{\alpha_{2},\alpha_1} \\
F_{\alpha_1,\alpha_{n-1}} & F_{\alpha_1,\alpha_{n-2}} & F_{\alpha_1,\alpha_{n-3}}  & \ldots & F_{\alpha_1,\alpha_{2}} & F_{\alpha_1,\alpha_{1}}  \\
\end{pNiceMatrix}$, so its characteristic polynomial is $p(\lambda)=(-1)^{q-\frac{n(n+1)}{2}}(q-1)^{n-1}q^\frac{(n-1)(n-2)}{2}\lambda^{q^{n-1}-n}\det(B-\lambda C)$, where $B \in M_{n-1}(\RR)$ is the upper-triangular matrix with all its elements equal to $1$ and $C=\antidiag(1/\alpha_1,1/\alpha_2,\ldots,1/\alpha_{n-1}) \in M_{n-1}(\RR)$ by the same proof as in the proof of Theorem \ref{xjenula}. 
Compare this result with \cite[Corollary 2.8 and Corollary 2.13]{pirzada}, where similar results (but for rings $\ZZ_{p^nq^m}$) have been achieved for zero-divisor graphs (without loops), using completely different methods.
\end{Remark}

We have the following structural lemma, the proof of which is implicit in the proof of \cite[Theorem 4.6]{dolzan1}, but we include it here for the sake of completeness.

\begin{Lemma}
    \label{structure}
    Let $R$ be a finite commutative local ring of order $q^n$, $R/J$ a field of order $q$ and $J^{n-1} \neq 0$. Then there exist $g \in R^*$ and $x \in J \setminus J^2$ such that $R=\{\sum_{i=0}^{n-1}\lambda_i x^i;  \lambda_i \in \{0, g, g^2, \ldots, g^{q-1} \} \text { for every } i \in \{0,1,\ldots,n-1\}\}$. Furthermore, if $g^k-g^l \in J$ for some $k,l \in \{0,1,\ldots,q-1\}$, then $k=l$.
\end{Lemma}
\begin{Proof}
    Choose $g_1 \in R^*$ such that $\overline{g_1}$ is a cyclic generator of the field
    $R/J$. Since $q=p^r$ for some prime $p$ and $|R^*|=q^{n-1}(q-1)$, the multiplicative order of $g_1$ is $p^{t}(q-1)$ for some integer $t$. Now, denote $g=g_1^t$. Since $q-1$ is prime to $p$, $g$ is an element of order $q-1$. Also, $\overline{g}$ is a cyclic generator of $R/J$. Therefore, if we have $g^k-g^l \in J$ for some $k,l \in \{0,1,\ldots,q-1\}$, we must have $k=l$.
    Now, choose $x_i \in J^i \setminus J^{i+1}$ for every $i \in \{0,1,\ldots,n-1\}$. Assume that 
    $\sum_{i=0}^{n-1}\lambda_ix_i=\sum_{i=0}^{n-1}\mu_ix_i$, where $\lambda_i, \mu_i \in \{0, g, g^2, \ldots, g^{q-1} \}$ for every $i \in \{0,1,\ldots,n-1\}$.
    Since $(\lambda_0-\mu_0)x_0 \in J$ for $x_0 \notin J$, we have $\lambda_0=\mu_0$. We then proceed by induction --- assume that $\lambda_i=\mu_i$ for $i=0,1,\ldots,j-1$. Then $\sum_{i=j}^{n-1}(\lambda_i-\mu_i)x_i=0$, so $(\lambda_{j}-\mu_{j})x_j \in J^{j+1}$ for $x_j \notin J^{j+1}$, so we have $\lambda_j=\mu_j$.
   This implies $\lambda_i=\mu_i$ for every $i \in \{0, 1,\ldots,n-1\}$, so $R=\{\sum_{i=0}^{n-1}\lambda_i x_i;  \lambda_i \in \{0, g, g^2, \ldots, g^{q-1} \} \text { for every } i \in \{0,1,\ldots,n-1\}\}$ and $J=\{\sum_{i=1}^{n-1}\lambda_i x_i;  \lambda_i \in \{0, g, g^2, \ldots, g^{q-1} \} \text { for every } i \in \{1,2,\ldots,n-1\}\}$. Since $J^{n-1} \neq 0$ and $J^n=0$, we see that $x_1^{n-1} \neq 0$, so we can choose $x=x_1$ and $x_i=x^i$ and the lemma is proven.
\end{Proof}

We now examine the case $u \in R^*$ when the order of $R$ is even.

\begin{Theorem}
  \label{xjeenotaeven}
   Let $R$ be a finite commutative local ring of order $q^n$ with $J^{n-1} \neq 0$ such that $R/J$ is a field of even order $q$. For $u \in R^*$, let $p(\lambda)$ denote the characteristic polynomial of the matrix $A_u(R)$. Then the following statements hold.
   \begin{enumerate}
       \item 
   If the characteristic of $R$ is $2$, then $p(\lambda)=\lambda^{q^{n-1}}(\lambda+1)^{q^n-q^{n-1}}$.
       \item 
   Suppose that the characteristic of $R$ is $2^n$ with $n \geq 3$. If $u=s^2$ for some $s \in R^*$, then  $p(\lambda)=\lambda^{q^{n-1}}(\lambda-1)^{\frac{q^n-q^{n-1}+4}{2}}(\lambda+1)^{\frac{q^n-q^{n-1}-4}{2}}$. For all other $u \in R^*$, we have $p(\lambda)=\lambda^{q^{n-1}}(\lambda-1)^{\frac{q^n-q^{n-1}}{2}}(\lambda+1)^{\frac{q^n-q^{n-1}}{2}}$.
   Furthermore, there exist exactly $\frac{q^n-q^{n-1}}{4}$ elements $u \in R^*$ such that $u=s^2$ for some $s \in R^*$.
\end{enumerate}
\end{Theorem}
\begin{Proof}
  We proceed similarly as in the proof of Theorem \ref{xjeenota}. Again, we have to consider the case when there exists an element $s \in R^*$ such that $s^{-1}u=s$, that is $s^2=u$, and the case when $u$ is not a square.
    
  Assume firstly that $\ch(R)=2$. If $s^2=v^2$ for some $s, v \in R^*$, then the fact that $R/J$ is a field of characteristic $2$ implies that $s=v+j$ for some $j \in J$. 
 Then $s^2=v^2+j^2=v^2$, so $j^2=0$. Obviously, if $s=v+j$ with $j^2=0$, we also have $s^2=v^2$.
 Lemma \ref{structure} tells us that there exist $g \in R^*$ and $x \in J$ such that $R=\{\sum_{i=0}^{n-1}\lambda_i x^i;  \lambda_i \in \{0, g, g^2, \ldots, g^{q-1} \} \text { for every } i \in \{0,1,\ldots,n-1\}\}$. Observe that since $J^{n-1} \neq 0$, we have $\{j \in J; j^2=0\}=\{\sum_{i=\lceil\frac{n}{2}\rceil}^{n-1}\lambda_i x^i;  \lambda_i \in \{0, g, g^2, \ldots, g^{q-1} \} \text { for every } i \in \{\lceil\frac{n}{2}\rceil,\lceil\frac{n}{2}\rceil+1,\ldots,n-1\}\}$. Thus,  we have proved that for any $v \in R$ we have $|\{s \in R;s^2=v^2\}|=|\{\sum_{i=\lceil\frac{n}{2}\rceil}^{n-1}\lambda_i x^i;  \lambda_i \in \{0, g, g^2, \ldots, g^{q-1} \} \text { for every } i \in \{\lceil\frac{n}{2}\rceil,\lceil\frac{n}{2}\rceil+1,\ldots,n-1\}\}|=q^{n-\lceil\frac{n}{2}\rceil}=q^{\lfloor\frac{n}{2}\rfloor}$. 
 Again, denote $B_k=\begin{pNiceMatrix}
    0 & 1 \\

1 &  0 \end{pNiceMatrix} \oplus \begin{pNiceMatrix}
    0 & 1 \\

1 &  0 \end{pNiceMatrix} \oplus \ldots \oplus \begin{pNiceMatrix}
    0 & 1 \\

1 &  0 \end{pNiceMatrix} \in M_{2k}(\RR)$ ($k$ summands).
   This means that in the case when $u$ is a square, we have (with the reordering of elements in $R$: $q^{n-1}$ elements from $J$, $q^{\lfloor \frac{n}{2}\rfloor}$ elements $s \in R \setminus J$ such that $s^2=u$, and finally the other $q^n-q^{n-1}-q^{\lfloor \frac{n}{2}\rfloor}$ invertible elements from $R$ in pairs $s$, $s^{-1}u$),
    $A_u(R)=0_{q^{n-1}} \oplus I_{q^{\lfloor\frac{n}{2}\rfloor}} \oplus B_\frac{q^n-q^{n-1}-q^{\lfloor\frac{n}{2}\rfloor}}{2}$. Thus, we have $p(\lambda)=\lambda^{q^{n-1}}(\lambda+1)^{q^n-q^{n-1}}$ 
    in this case. 
        In the case when $u$ is not a square, we have (with a similar reordering),
$A_u(R)=0_{q^{n-1}} \oplus B_\frac{q^n-q^{n-1}}{2}$, therefore again $p(\lambda)=\lambda^{q^{n-1}}(\lambda+1)^{q^n-q^{n-1}}$.

  Assume now that $\ch(R)=2^n$ and $n > 1$. Then $2 \in J \setminus J^2$, so by Lemma \ref{structure}, we have $g \in R^*$ such that
  $R=\{\sum_{i=0}^{n-1}\lambda_i 2^i;  \lambda_i \in \{0, g, g^2, \ldots, g^{q-1} \} \text { for every } i \in \{0,1,\ldots,n-1\}\}$ and $J=(2)$.

  
  Choose $s=\sum_{i=0}^{n-1}\lambda_i 2^i$ and observe that $s^2=\sum_{i=0}^{\lfloor\frac{n-1}{2}\rfloor}\lambda_i^2 2^{2i}+\sum_{i=0}^{n-2}\sum_{j=0}^{\lfloor\frac{i}{2}\rfloor}\lambda_j\lambda_{i-j}2^{i+1}$. 
  Suppose now that $s^2=v^2$ for some $v=\sum_{i=0}^{n-1}\lambda_i' 2^i \in R^*$.
  Then $\lambda_0^2=\lambda_0'^2$, so $\lambda_0'=\pm \lambda_0 \neq 0$. Note that $2\lambda_0\in J$, so $-\lambda_0 \in \lambda_0 + J$, hence we can assume that $\lambda_0'=\lambda_0$. The next equation (when comparing coefficients at $2^2$) is $\lambda_0\lambda_1+\lambda_1^2=\lambda_0'\lambda_1'+\lambda_1'^2$. Since $\lambda_0=\lambda_0'$, we are looking for all the distinct roots of the polynomial $f(x)=x^2+\lambda_0x-\lambda_0\lambda_1-\lambda_1^2 \in R[x]$, but again only modulo $J$. Obviously, $x=\lambda_1$ and $x=-\lambda_1-\lambda_0$ are two distinct roots (since $0 \neq \lambda_0=-2\lambda_1 \in J$ is not possible), so we have two possibilities for choosing $\lambda_1'$. 
  
  Suppose now that $\lambda_0',\lambda_1',\ldots,\lambda_{t-1}'$ are already determined. Then by observing the coefficient at $2^{t+1}$ (for $t \leq n-2$) in the equation $s^2=v^2$, we have $\sum_{j=0}^{\lfloor\frac{t}{2}\rfloor}\lambda_j\lambda_{t-j}=\sum_{j=0}^{\lfloor\frac{t}{2}\rfloor}\lambda_j'\lambda_{t-j}'$ for even numbers $t$ and $\lambda_{\frac{t+1}{2}}^2+\sum_{j=0}^{\lfloor\frac{t}{2}\rfloor}\lambda_j\lambda_{t-j}=\lambda_{\frac{t+1}{2}}'^2+\sum_{j=0}^{\lfloor\frac{t}{2}\rfloor}\lambda_j'\lambda_{t-j}'$ for odd numbers $t$. Note that both of these two equations uniquely determine the value of $\lambda_{t+1}'$ for any $2 \leq t \leq n-2$. Obviously, $\lambda'_{n-1}$ can be arbitrary since it disappears when squared. We have therefore proved that for any $s \in R^*$ there exist exactly $4$ elements $v \in R^*$ such that $s^2=v^2$. This further implies that the mapping $\alpha : R^* \rightarrow R^*$ defined by $\alpha(x)=x^2$ has an image of order $\frac{|R^*|}{4}=\frac{q^n-q^{n-1}}{4}$.
   So, if $u$ is a square, we have (again, with the reordering of elements in $R$: $q^{n-1}$ elements from $J$, $4$ elements $s \in R \setminus J$ such that $s^2=u$, and finally the other $q^n-q^{n-1}-4$ invertible elements from $R$ in pairs $s$, $s^{-1}u$),
    $A_u(R)=0_{q^{n-1}} \oplus I_{4} \oplus B_\frac{q^n-q^{n-1}-4}{2}$ and thus $p(\lambda)=\lambda^{q^{n-1}}(\lambda-1)^{\frac{q^n-q^{n-1}+4}{2}}(\lambda+1)^{\frac{q^n-q^{n-1}-4}{2}}$. In the case when $u$ is not a square, we have (with a suitable reordering of elements),
$A_u(R)=0_{q^{n-1}} \oplus B_\frac{q^n-q^{n-1}}{2}$, therefore $p(\lambda)=\lambda^{q^{n-1}}(\lambda-1)^{\frac{q^n-q^{n-1}}{2}}(\lambda+1)^{\frac{q^n-q^{n-1}}{2}}$.
  \end{Proof}

\begin{Example}
    Let us illustrate Theorem \ref{xjeenotaeven} by two examples.
    \begin{enumerate}
        \item 
        Let $n$ be an integer. Let $T$ denote the $n$-by-$n$ matrix with $1$ on the first upper diagonal and $0$ elsewhere.
        Then observe that
        $R=\{a_0I+a_1T+a_2T^2+\ldots+a_{n-1}T^{n-1}; a_0,a_1,\ldots,a_{n-1} \in \ZZ_2\} \subseteq M_n(\ZZ_{2})$ is a finite local commutative ring of order $2^n$ with $J(R)^{n-1} \neq 0$. Since the characteristic of $R$ is $2$, we have by Theorem \ref{xjeenotaeven}(1) that $p(\lambda)=\lambda^{2^{n-1}}(\lambda+1)^{2^n-2^{n-1}}$ for any $u \in R^*$.
        \item
        If $R=\ZZ_{8}$, then Theorem \ref{xjeenotaeven}(2) gives us 
        $p(\lambda)=\lambda^4(\lambda-1)^4$ for $u=1$ and
        $p(\lambda)=\lambda^4(\lambda-1)^2(\lambda+1)^2$ for $u=3,5,7$.
    \end{enumerate}
\end{Example}

\bigskip

%
%

The following lemma is common knowledge (see, for example, \cite[p. 5]{schur}).

\begin{Lemma}(Schur's formula)
    \label{detident} Let $A \in M_n(\RR), D \in M_m(\RR)$ and let $B, C$ be matrices of appropriate sizes such that
    $M=\begin{pNiceMatrix}
    A & B \\

C &  D \end{pNiceMatrix} \in M_{n+m}(\RR)$. If $A$ is an invertible matrix, then $\det(M)=\det(A)\det(D-CA^{-1}B)$.
\end{Lemma}

Finally, if $u \in Z(R) \setminus \{0\}$, we have the following theorem.

\begin{Theorem}
  \label{xjedn}
   Let $R$ be a finite commutative local ring of order $q^n$, $R/J$ a field of order $q$ and $J^{n-1} \neq 0$. Let $k \geq 0$ be such that $u \in J^k \setminus J^{k+1}$ and let $p(\lambda)$ denote the characteristic polynomial of the matrix $A_u(R)$. If $k$ is odd, then $p(\lambda)=(-1)^q\lambda^{q^{n-k-1}(q^{k+1}-(k+1)(q-1))}(\lambda^2-q^{k})^{q^{n-k-1}(q-1)\frac{k+1}{2}}$. If $k$ is even and $q$ is odd, then
   $p(\lambda)=-\lambda^{q^{n-k-1}(q^{k+1}-(k+1)(q-1))}(\lambda^2-q^{k})^{\frac{q^{n-k-1}(q-1)}{2}(k+1)-1}(\lambda - q^{k/2})^2$ if there exists $s \in R$ such that $u=s^2$, and $p(\lambda)=-\lambda^{q^{n-k-1}(q^{k+1}-(k+1)(q-1))}(\lambda^2-q^{k})^{q^{n-k-1}(q-1)\frac{k+1}{2}}$ otherwise.
  Furthermore, there exist exactly $\frac{q^{n-k-1}(q-1)}{2}$ elements $u \in J^k \setminus J^{k+1}$, such that $u=s^2$ for some $s \in R$.
\end{Theorem}
\begin{Proof}
     By Lemma \ref{structure}, there exist $g \in R^*$ and $x \in J$ such that $R=\{\sum_{i=0}^{n-1}\lambda_i x^i;  \lambda_i \in \{0, g, g^2, \ldots, g^{q-1} \} \text { for every } i \in \{0,1,\ldots,n-1\}\}$. Therefore $u=x^k \sum_{i=0}^{n-k-1}\lambda_{k+i} x^i$ where $\lambda_k \neq 0$. Obviously, for $s \in J^{l}$ with $l \geq k+1$ there exists no $t \in R$ such that $st=u$, so we assume from now on that $l \leq k$. In fact, since $st=ts$ and $s \in J^l \setminus J^{l+1}$ implies that $t \in J^{k-l} \setminus J^{k-l+1}$, we can further limit ourselves to the cases when $l \leq k/2$.
          
     So, suppose that $s \in J^l \setminus J^{l+1}$ for some integer $l \leq k/2$. We have $s=x^l \sum_{i=0}^{n-l-1}\mu_{l+i} x^i$ where $\mu_l \neq 0$. Denote $s'=\sum_{i=0}^{n-l-1}\mu_{l+i} x^i$ and $u'=\sum_{i=0}^{n-k-1}\lambda_{k+i} x^i$ and observe that $s', u' \in R^*$. Thus, $st=u$  implies that $x^l t=x^ku's'^{-1}$, so $t \in J^{k-l} \setminus J^{k-l+1}$. Write $t=x^{k-l}t'$ for some $t' \in R^*$, which yields 
     \begin{equation}
     \label{eq1}
         x^k(t'-u's'^{-1})=0.
     \end{equation} 
     Denote $T'(s)=\{t \in J^{k-l} \setminus J^{k-l+1}; st=u\}$. By the above, $t'=u's'^{-1}+t''$ for any $t'' \in J^{n-k}$ are all the solutions of equation (\ref{eq1}), therefore $T'(s)=\{x^{k-l}u's'^{-1}+t'''; t''' \in J^{n-l}\}$. Note also that $|T'(s)|=q^{l}$.

     Choose $s_1=x^ls_1', s_2=x^ls_2' \in J^l \setminus J^{l+1}$. If $s_1-s_2 \in J^{n-k+l}$, then $s_1'-s_2' \in J^{n-k}$ and since $J^{n-k}$ is an ideal in $R$, we have $s_1'^{-1}-s_2'^{-1} \in J^{n-k}$, so $s_2'^{-1}=s_1'^{-1}+t$ for some $t \in J^{n-k}$. This implies that $T'(s_2)=\{x^{k-l}u's_2'^{-1}+t'''; t''' \in J^{n-l}\}=\{x^{k-l}u's_1^{-1}+t'''; t''' \in J^{n-l}\}=T'(s_1)$.
     
     Now, suppose that $t \in T'(s_1) \cap T'(s_2)$. Then $t=x^{k-l}u's_1'^{-1}+t_1$ for some $t_1 \in J^{n-l}$ and $t=x^{k-l}u's_2'^{-1}+t_2$ for some $t_2 \in J^{n-l}$. Thus, $x^{k-l}(u's_1'^{-1}-u's_2'^{-1}) \in J^{n-l}$, so $s_1'^{-1}-s_2'^{-1} \in J^{n-k}$ and therefore $s_1'-s_2' \in J^{n-k}$, which implies that $s_1-s_2 \in J^{n-k+l}$, and thus $T'(s_1)=T'(s_2)$. 
     
     We have therefore proved that $T'(s_1)=T'(s_2)$ if and only if $s_1-s_2 \in J^{n-k+l}$. Since $|J^{n-k+l}|=q^{k-l}$, all the $q^{n-l-1}(q-1)$ elements in $J^l \setminus J^{l+1}$ are divided into disjoint sets $S_1,S_2, \ldots, S_d$ of orders $q^{k-l}$ (where $d=\frac{q^{n-l-1}(q-1)}{q^{k-l}}=q^{n-k-1}(q-1)$), such that $T'(s)=T'(s')$ if and only if $s, s' \in S_j$ for some $j \in \{1,2,\ldots,d\}$.
          
     If $k$ is odd, then $k \neq 2l$ for all $l$, so $s^2 \neq u$ for any $s \in R$. Therefore, $s \notin T'(s)$, so (with the following reordering of elements in $R$: $q^{n-k-1}$ elements from $J^{k+1}$, then for every $l=0$ to $\frac{k-1}{2}$, we have first the elements from $S_i \subseteq J^l \setminus J^{l+1}$ and then the elements from $T'(s) \subseteq J^{k-l} \setminus J^{k-l+1}$ for some $s \in S_i$, for all $i=1,2,\ldots,q^{n-k-1}(q-1)$),
     we get $A_u(R)=0_{q^{n-k-1}} \oplus_{l=0}^{\frac{k-1}{2}} \oplus_{r=1}^{q^{n-k-1}(q-1)} C_{q^{k-l},q^l}$, where for any integers $i,j$, $C_{i,j}$ is a block matrix $\begin{pNiceMatrix}
    0_{i} & F_{i,j} \\

F_{j,i} &  0_{j} \end{pNiceMatrix}$. Lemma \ref{detident} now yields $\det(C_{i,j}-\lambda I)=\det(\lambda^2I-F_{i,j}F_{j,i})$. The fact that $F_{i,j}F_{j,i}=jF_{i,i}$ which has the only non-zero eigenvalue $ij$, now ensures us that the characteristic polynomial of $C_{q^{k-l},q^l}$ is $(\lambda^2-q^k)\lambda^{q^{k-l}+q^l-2}$. Thus, the characteristic polynomial of the matrix $A_u(R)$ is $p(\lambda)=(-1)^q\lambda^{q^{n-k-1}(q^{k+1}-(k+1)(q-1))}(\lambda^2-q^{k})^{q^{n-k-1}(q-1)\frac{k+1}{2}}$.

Suppose now that $k$ is even and $q$ is odd. If $k \neq 2l$, then $s^2 \neq u$ for any $s \in J^l \setminus J^{l+1}$, so the same reasoning as above still applies. Consider now the case $k=2l$.  
Let $s=x^ls'$ for some $s' \in R^*$ and suppose that $s^2=u$. 
This implies that $(s'^2-u')x^k=0$, so $u'=s'^2+t$ for some $t \in J^{n-k}$. Suppose that $s^2=v^2$ for some $v \in R$. Therefore $v=x^lv'$ for some $v' \in R^*$, so $s'^2+t=v'^2+t'$ for some $t' \in J^{n-k}$. Since $R/J$ is a field, we have $s'=\pm v'+j$ for some $j \in J$.    But then $s'^2=v'^2 \pm 2v'j + j^2=v'^2 + j(\pm 2v' + j)$. Since $q$ is odd, $2 \in R^*$, so $\pm 2v' + j \in R^*$ and therefore $j \in  J^{n-k}$.  This implies that $s=\pm v+j'$ for some $j' \in J^{n+l-k}=J^{n-l}$. But now $s^2=v^2$ implies $v^2 \pm 2vj' + j'^2=v^2$. Since $v \in J^l$ and $j' \in J^{n-l}$, we have $vj'=0$, therefore $j'^2=0$. We have proved that $j' \in J^{\max\{n-l,\lceil n/2 \rceil\}}$. However $k=2l < n$, so $l \leq \lfloor n/2 \rfloor$ and consequently $n-l \geq \lceil n/2 \rceil$. Therefore, we have proved that $s=\pm v + j'$ for some $j' \in J^{n-l}$. Conversely, for every $j' \in J^{n-l}$ we have $vj'=0$, and by the fact that $n-l \geq \lceil n/2 \rceil$ also $j'^2=0$, which has a consequence that $(\pm v+j')^2=v^2$. 
So, if $u=s^2$ for some $s \in R$, there exist exactly $2q^l$ elements $v \in J^l$ such that $u=v^2$. Note that if $s^2=u$ and $st=u$, then $s \in T'(s) \cap T'(t)$, which by the above implies that $T'(s) = T'(t)$, so $t \in T'(s)$ implies $t^2=u$ as well. Since $q$ is odd, $2$ is invertible and thus $v^2 \neq -v^2$, so this implies that these $2q^l$ elements $v$ such that $u=v^2$ are divided into exactly two sets $S_1, S_2$ of order $q^l$ such that for $i \in \{1,2\}$ and 
for any $x \in S_i$ we have $xy=u$ if and only if $y \in S_i$.
Furthermore, we are left with exactly $q^{n-l-1}(q-1)-2q^l$ elements $s \in J^l \setminus J^{l+1}$ such that $s^2 \neq u$. By the same reasoning as above, these elements are divided into $d=\frac{q^{n-l-1}(q-1)-2q^l}{q^{k-l}}=q^{n-k-1}(q-1)-2$ disjoint sets $S_1,S_2, \ldots, S_d$ of orders $q^l$ such that 
$T'(s)=T'(s')$ if and only if $s, s' \in S_j$ for some $j \in \{1,2,\ldots,d\}$. Since for every such $s$, the set $T'(s)$ is in this case a subset of $J^l \setminus J^{l+1}$, the $\frac{q^{n-k-1}(q-1)-2}{2}$ different pairs of these sets yield  the product matrices $C_{q^l,q^l}$.
Therefore, if $u \in J^k \setminus J^{k+1}$ is a square, we reorder the elements in $R$ thus: as above, we start with the $q^{n-k-1}$ elements from $J^{k+1}$, then for every $l=0$ to $\frac{k}{2}-1$, we have first the elements from $S_i \subseteq J^l \setminus J^{l+1}$ and then the elements from $T'(s) \subseteq J^{k-l} \setminus J^{k-l+1}$ for some $s \in S_i$, for all $i=1,2,\ldots,q^{n-k-1}(q-1)$.  But this time we are still left with all the elements in $J^{\frac{k}{2}} \setminus J^{\frac{k}{2}+1}$. By the above reasoning, we can divide them first into two sets of cardinality $q^{\frac{k}{2}}$ such that all the products of elements from each set yield $u$, and then the rest of the elements are divided into sets $S_1, S_2, \ldots, S_{q^{n-k-1}(q-1)-2}$ of cardinality $q^{\frac{k}{2}}$, but by the preceding arguments, we can order them in such a way that the products $S_i S_{i+1}$ yield $u$ for all odd numbers $i$). Therefore, we get $A_u(R)=0_{q^{n-k-1}} \oplus_{l=0}^{\frac{k}{2}-1} \oplus_{r=1}^{q^{n-k-1}(q-1)} C_{q^{k-l},q^l} \oplus F_{q^{k/2}} \oplus F_{q^{k/2}} \oplus_{r=1}^{\frac{q^{n-k-1}(q-1)}{2}-1} C_{q^{\frac{k}{2}},q^{\frac{k}{2}}}$. Since the only non-zero eigenvalue of $F_{q^{k/2}}$ is $q^{k/2}$, this implies that
$p(\lambda)=-\lambda^{q^{n-k-1}(q^{k+1}-(k+1)(q-1))}(\lambda^2-q^{k})^{\frac{q^{n-k-1}(q-1)}{2}(k+1)-1}(\lambda - q^{k/2})^2$.

Finally, in the case that $u$ is not a square, we have in the case $k=2l$ that there are exactly $q^{n-l-1}(q-1)$ elements $s \in J^l \setminus J^{l+1}$, divided into $d=\frac{q^{n-l-1}(q-1)}{q^{k-l}}=q^{n-k-1}(q-1)$ disjoint sets $S_1,S_2,\ldots,S_d$ of orders $q^l$ such that 
$T'(s)=T'(s')$ if and only if $s, s' \in S_j$ for some $j \in \{1,2,\ldots,d\}$. Since again for every $s \in J^l \setminus J^{l+1}$, the set $T'(s)$ is a subset of $J^l \setminus J^{l+1}$, the $\frac{q^{n-k-1}(q-1)}{2}$ different pairs of these sets yield the product matrices $C_{q^l,q^l}$. Therefore (by a similar reordering of $R$ as before), $A_u(R)=0_{q^{n-k-1}} \oplus_{l=0}^{\frac{k}{2}-1} \oplus_{r=1}^{q^{n-k-1}(q-1)} C_{q^{k-l},q^l}\oplus_{r=1}^{\frac{q^{n-k-1}(q-1)}{2}} C_{q^{\frac{k}{2}},q^{\frac{k}{2}}}$. In this case, we then have
$p(\lambda)=-\lambda^{q^{n-k-1}(q^{k+1}-(k+1)(q-1))}(\lambda^2-q^{k})^{q^{n-k-1}(q-1)\frac{k+1}{2}}$.

Note further that by above, in the case $k=2l$ the 
mapping $\alpha: J^l \setminus J^{l+1} \mapsto J^k \setminus J^{k+1}$ defined by $\alpha(x)=x^2$ has an image of order $\frac{q^{n-l-1}(q-1)}{2q^l}=\frac{q^{n-k-1}(q-1)}{2}$, so exactly half of the elements in $J^k \setminus J^{k+1}$ are squares.
     
    
    
\end{Proof}

\begin{Example}
    Let us again illustrate Theorem \ref{xjedn} by two examples. Choose a prime $q$ and integers $0 \leq k \leq n$. Let $R=\ZZ_{q^n}$, $u=\alpha q^k$ for some integer $\alpha$ that is not divisible by $q$. Let $p(\lambda)$ denote the characteristic polynomial of the matrix $A_u(R)$.
    \begin{enumerate}
        \item 
        If $k$ is odd, then we have $p(\lambda)=(-1)^q\lambda^{q^{n-k-1}(q^{k+1}-(k+1)(q-1))}(\lambda^2-q^{k})^{q^{n-k-1}(q-1)\frac{k+1}{2}}$.
        \item
        Suppose $q=3$, $k=2$ and $n=3$.
        Then $p(\lambda)=-\lambda^{21}(\lambda+3)^{2}(\lambda - 3)^4$ if $u=9$, and $p(\lambda)=-\lambda^{21}(\lambda+3)^3(\lambda-3)^3$ if $u=18$.
    \end{enumerate}
\end{Example}

\bigskip

 \section{Local rings with the Jacobson radical of minimal nilpotency index}

\bigskip

In this section, we will concern ourselves with finite commutative local rings $R$ such that $J^{2} = 0$. Our main result is the following.

\begin{Theorem}
  \label{jkvadratnula}
   Let $R$ be a finite commutative local ring of order $q^n$, $R/J$ a field of order $q$ and $J^{2} = 0$. Let $p(\lambda)$ denote the characteristic polynomial of the matrix $A_u(R)$. Then $p(\lambda)=
   \begin{cases}
   (-1)^q\lambda^{q^n-3}(\lambda^3-q^{n-1}\lambda^2-(q^n-q^{n-1})\lambda+q^{n-1}(q-1)(q^{n-1}-1)), \text{ if } u=0; \\
   \lambda^{q^{n-1}}(\lambda-1)^\frac{q^n}{2}(\lambda+1)^{\frac{q^n-2q^{n-1}}{2}}, \text{ if } q \text { is even and } u = s^2 \text{ for some } s \in R^*; \\
   \lambda^{q^{n-1}}(\lambda-1)^\frac{q^n-q^{n-1}}{2}(\lambda+1)^{\frac{q^n-q^{n-1}}{2}}, \text{ if } q \text { is even}, u \in R^* \text { and } u \neq s^2 \text{ for all } s \in R^*; \\
   -\lambda^{q^{n-1}}(\lambda-1)^{\frac{q^n-q^{n-1}+2}{2}}(\lambda+1)^{\frac{q^n-q^{n-1}-2}{2}} \text{ if } q \text { is odd and } u = s^2 \text{ for some } s \in R^*; \\
   -\lambda^{q^{n-1}}(\lambda-1)^{\frac{q^n-q^{n-1}}{2}}(\lambda+1)^{\frac{q^n-q^{n-1}}{2}}, \text{ if } q \text { is odd}, u \in R^* \text { and } u \neq s^2 \text{ for all } s \in R^*; \\
   (-1)^q\lambda^{q^n-2q+2}(\lambda^2-q^{n-1})^{q-1}, \text{ if } 0 \neq u \in J.
   \end{cases}$
\end{Theorem}
\begin{Proof}
  Examine firstly the case $u=0$. If $s \in R^*$, then $st=u$ if and only if $t=0$. If $s \in J \setminus \{0\}$, then $st=u$ if and only if $t \in J$, and finally if $s = 0$, then obviously $sR=u$. We reorder the elements in $R$: $R^*$, $J\setminus \{0\}, 0$. This gives us that $A_u(R)$ is equal to the block matrix $\begin{pNiceMatrix}
    0_{q^n-q^{n-1}} & 0_{q^n-q^{n-1},q^{n-1}-1} & F_{q^n-q^{n-1},1} \\
    0_{q^{n-1}-1, q^n-q^{n-1}} & F_{q^{n-1}-1, q^{n-1}-1} & F_{q^{n-1}-1,1} \\
   F_{1, q^n-q^{n-1}} & F_{1, q^{n-1}-1} & F_{1,1}
\end{pNiceMatrix}$. Thus $0$ is an eigenvalue of $A_u(R)$ with multiplicity $q^n-3$. Define now the vectors $v_0=(0_{1,q^n-q^{n-1}},0_{1, q^{n-1}-1},1)^T, v_1=(0_{1,q^n-q^{n-1}},F_{1, q^{n-1}-1},1)^T$ and $v_2=(F_{1,q^n-q^{n-1}},F_{1, q^{n-1}-1},1)^T$. 
Observe that $A_u(R)v_0=v_2$, $A_u(R)v_1=(q^{n-1}-1)v_1+v_2$ and $A_u(R)v_2=(q^n-q^{n-1})v_0+(q^{n-1}-1)v_1+v_2$. Choose $w=\alpha_0v_0+\alpha_1v_1+\alpha_2v_2$ for some $\alpha_0, \alpha_1, \alpha_2 \in \CC$ and observe that $A_u(R)w=\lambda w$ in the basis $\{v_0,v_1,v_2\}$ transforms into  
$\begin{pNiceMatrix}
    0 & 0 & q^n-q^{n-1} \\
    0 & q^{n-1}-1 & q^{n-1}-1 \\
    1 & 1 & 1
\end{pNiceMatrix}
\begin{pNiceMatrix}
    \alpha_0 \\
    \alpha_1 \\
   \alpha_2
\end{pNiceMatrix}=\lambda\begin{pNiceMatrix}
    \alpha_0 \\
    \alpha_1 \\
   \alpha_2
\end{pNiceMatrix}$. Calculating the characteristic polynomial of this $3$-by-$3$ matrix, we get $-\lambda^3+q^{n-1}\lambda^2+(q^n-q^{n-1})\lambda-q^{n-1}(q-1)(q^{n-1}-1)$. Therefore $p(\lambda)=(-1)^q\lambda^{q^n-3}(\lambda^3-q^{n-1}\lambda^2-(q^n-q^{n-1})\lambda+q^{n-1}(q-1)(q^{n-1}-1))$ is the characteristic polynomial of $A_u(R)$ in this case.

Assume now that $u \in R^*$ and $q$ is even. Obviously, for $s \in J$ there exists no $t \in R$ such that $st=u$. If $s \in R^*$, then $t=s^{-1}u$ is the only solution to the equation $st=u$. Since $2 \in J$ and $J^2=0$, we have $(x+j)^2=x^2$ for any $x \in R$ and any $j \in J$. Furthermore, if $x^2=y^2$ for some $x,y \in R^*$, then the fact that $R/J$ is a field of characteristic $2$ implies that $y=x + j$ for some $j \in J$.  Thus, if $u$ is a square, we have (with the reordering of elements in $R$: $J$, then the $q^{n-1}$ elements $s \in R^*$ such that $s^2=u$, and finally the other $q^n-2q^{n-1}$ elements from $R^*$ in pairs $s$, $s^{-1}u$), $A_u(R)=0_{q^{n-1}} \oplus I_{q^{n-1}} \oplus B_\frac{q^n-2q^{n-1}}{2}$, therefore $p(\lambda)=\lambda^{q^{n-1}}(\lambda-1)^\frac{q^n}{2}(\lambda+1)^{\frac{q^n-2q^{n-1}}{2}}$. And, if $u$ is not a square, then (with a similar reordering as above), $A_u(R)=0_{q^{n-1}} \oplus B_\frac{q^n-q^{n-1}}{2}$, so $p(\lambda)=\lambda^{q^{n-1}}(\lambda-1)^\frac{q^n-q^{n-1}}{2}(\lambda+1)^{\frac{q^n-q^{n-1}}{2}}$.

If $u \in R^*$ and $q$ is odd, then the result follows by Theorem \ref{xjeenota}.

Finally, assume that $0 \neq u \in J$. If $s \in R^*$, then $st=u$ has a unique solution $t=s^{-1}u$ and obviously $s^{-1}u \neq s$. Note that $s^{-1}u=s'^{-1}u$ for some $s, s' \in R^*$ if and only if $(s^{-1}-s'^{-1})u=0$, which is true if and only if $s^{-1}-s'^{-1} \in J$. Since $J$ is an ideal, this holds if and only if $s-s' \in J$.
If $0 \neq s \in J$, the fact that $J^2=0$ shows that $st=u$ implies that $t \in R^*$. Furthermore, if $st=u$, then also $s(t+j)=u$ for any $j \in J$.
This time, we reorder the elements in $R$ thus: first all the invertible elements in $R$, divided into all non-zero cosets of the additive group $R/J$ --- $s_1+J, s_2+J, \ldots, s_{q-1}+J$, then the elements $s_1^{-1}u,s_2^{-1}u,\ldots,s_{q-1}^{-1}u$, and finally the other $q^{n-1}-q+1$ elements in any order. This implies that $A_u(R)$ is equal to the block matrix $\begin{pNiceMatrix}
    0_{q^n - q^{n-1}} & I_{q-1} \otimes v & 0_{q^n - q^{n-1},q^{n-1}-q+1}\\

I_{q-1} \otimes v^T &  0_{q-1} & 0_{q-1,q^{n-1}-q+1} \\
0_{q^{n-1}-q+1, q^n- q^{n-1}} &  0_{q^{n-1}-q+1, q-1} & 0_{q^{n-1}-q+1}
\end{pNiceMatrix}$, where $v=(1,1,\ldots,1)^T \in \RR^{q^{n-1}}$. Obviously, $A_u(R)$ is a matrix of rank $2q-2$, so there exist exactly $2q-2$ non-zero eigenvalues. 

Any non-zero eigenvalue $\lambda$ satisfies the equation
$\det\begin{pNiceMatrix}
    -\lambda I_{q^n - q^{n-1}} & I_{q-1} \otimes v \\

I_{q-1} \otimes v^T &  -\lambda I_{q-1}  \\
\end{pNiceMatrix}=0$. By Lemma \ref{detident}, this is equivalent to  $\det(-\lambda I_{q^n - q^{n-1}})\det((I_{q-1} \otimes v^T)(-\lambda^{-1} I_{q^n - q^{n-1}})(I_{q-1} \otimes v)+\lambda I_{q-1})$=0. This implies that $\det(v^TvI_{q-1}-\lambda^2 I_{q-1})=0$, therefore $\lambda^2=v^Tv=q^{n-1}$. So, $p(\lambda)=(-1)^q\lambda^{q^n-2q+2}(\lambda^2-q^{n-1})^{q-1}$.
\end{Proof}

\begin{Example}
        Let $n \geq 2$ be an integer and for any $1 \leq i,j \leq n$ let $E_{i,j}$ denote the $n$-by-$n$ matrix with $1$ at entry $(i,j)$ and $0$ elsewhere. Observe that
        $R=\{a_0I+a_1E_{1,2}+a_2E_{1,3}+\ldots+a_{n-1}E_{1,n-1}; a_0,a_1,\ldots,a_{n-1} \in \ZZ_2\} \subseteq M_n(\ZZ_{2})$ is a finite local commutative ring of order $2^n$ with $J(R)^2= 0$. For $u \in R$, let $p(\lambda)$ denote the characteristic polynomial of the matrix $A_u(R)$. 
        
        Then Theorem \ref{jkvadratnula} gives us
        $p(\lambda)=\lambda^{2^n-3}(\lambda^3- 2^{n-1} \lambda^2-2^{n-1}\lambda+2^{n-1}(2^{n-1}-1))$ for $u=0$, 
   $p(\lambda)=\lambda^{2^{n-1}}(\lambda-1)^{2^{n-1}}$ for $u=I$, and
   $p(\lambda)=\lambda^{2^{n-1}}(\lambda-1)^{2^{n-2}}(\lambda+1)^{2^{n-2}}$ for all $I \neq u \in R^*$.
\end{Example}

\bigskip

\bibliographystyle{amsplain}
\bibliography{biblio}

\bigskip

\end{document}